\pdfoutput=1
\RequirePackage{silence}
\documentclass[a4paper,11pt]{amsart}
\usepackage[hmarginratio={1:1},vmarginratio={1:1},lmargin=60.0pt,tmargin=60.0pt]{geometry}

\usepackage[numbers]{natbib}

\allowdisplaybreaks

\usepackage[utf8]{inputenc}

\usepackage{amssymb,amsmath,amsthm,amsfonts,mathrsfs,bbm}
\usepackage[table]{xcolor}
\usepackage{graphicx}
\usepackage{mathtools}

\usepackage{ytableau}

\definecolor{spinach}{RGB}{46,139,87}
\definecolor{tomato}{RGB}{255,99,71}
\definecolor{pumpkin}{RGB}{224,180,80}

\definecolor{orchid}{RGB}{143,40,194}
\definecolor{lava}{RGB}{207,16,32}
\definecolor{mydarkblue}{RGB}{10,10,150}

\definecolor{phantom}{RGB}{255,105,180}

\usepackage[yyyymmdd,hhmmss]{datetime}
\usepackage{todonotes}

\usepackage{xparse}

\usepackage{enumitem}
\setlist[enumerate]{itemsep=0.15cm,label=\emph{\upshape(\alph*)}}
\setlist[enumerate,2]{itemsep=0.15cm,label=\emph{\upshape(\roman*)}}

\let\emph\relax
\DeclareTextFontCommand{\emph}{\em}

\usepackage{array}
\newcolumntype{C}{>{$}c<{$}}

\DeclarePairedDelimiterX{\set}[1]{\{}{\}}{\setargs{#1}}
\NewDocumentCommand{\setargs}{>{\SplitArgument{1}{|}}m}{\setargsaux#1}
\NewDocumentCommand{\setargsaux}{mm}
{\IfNoValueTF{#2}{#1} {#1\,\delimsize|\,\mathopen{}#2}}

\usepackage[all]{xy}
\usepackage{tikz}
\usetikzlibrary{cd}
\usetikzlibrary{decorations}
\usetikzlibrary{decorations.markings}
\usetikzlibrary{decorations.pathreplacing}
\usetikzlibrary{decorations.pathmorphing}
\usetikzlibrary{arrows.meta,shapes,positioning,matrix,calc}
\usetikzlibrary{shapes.callouts}
\tikzset{anchorbase/.style={baseline={([yshift=-0.5ex]current bounding box.center)}},
tinynodes/.style={font=\tiny,text height=0.25ex,text depth=0.05ex},
smallnodes/.style={font=\scriptsize,text height=0.75ex,text depth=0.15ex},
crossline/.style={preaction={draw=white,line width=10.0pt,-},preaction={draw=black,line width=2.0pt,-}},
usual/.style={line width=2.0,color=black},
phantom/.style={line width=2.0,color=phantom,densely dotted},
crosslinep/.style={preaction={draw=white,line width=10.0pt,-},preaction={draw=phantom,densely dotted,line width=2.0pt,-}},
}
\tikzstyle directed=[postaction={decorate,decoration={markings,mark=at position #1 with {\arrow[line width=0.5mm]{>}}}}]
\tikzstyle rdirected=[postaction={decorate,decoration={markings,mark=at position #1 with {\arrow[line width=0.5mm]{<}}}}]

\newcommand{\ie}{\text{i.e.}}
\newcommand{\eg}{\text{e.g.}}
\newcommand{\cf}{\text{cf.}}

\newcommand{\muta}{\text{mutatis mutandis}}

\newcommand{\acts}{\centerdot}
\renewcommand{\dots}{\text{...}}

\newcommand{\placeholder}{{}_{-}}

\newcommand{\vcirc}{\circ_{v}}
\newcommand{\hcirc}{\circ_{h}}
\newcommand{\munit}{\mathbbm{1}}

\newcommand{\setstuff}[1]{\mathrm{#1}}
\newcommand{\catstuff}[1]{\mathbf{#1}}

\newcommand{\morstuff}[1]{\mathrm{#1}}

\newcommand{\idmor}[1]{1_{#1}}

\newcommand{\Aut}{\setstuff{Aut}}
\newcommand{\End}{\setstuff{End}}
\newcommand{\Hom}{\setstuff{Hom}}

\newcommand{\C}{\mathbb{C}}
\newcommand{\Z}{\mathbb{Z}}
\newcommand{\R}{\mathbb{R}}

\newcommand{\N}{\mathbb{Z}_{\geq 0}}

\newcommand{\qpar}{q}
\newcommand{\Zq}{\Z[\qpar,\qpar^{-1}]}
\newcommand{\Nq}{\N[\qpar,\qpar^{-1}]}
\newcommand{\Nqq}{1+\qpar\N[\qpar]}
\newcommand{\qnum}[1]{[#1]_{\qpar}}
\newcommand{\qfac}[1]{[#1]_{\qpar}!}
\newcommand{\qbin}[2]{{\textstyle\genfrac{[}{]}{0pt}{}{#1}{#2}}_{\qpar}}
\newcommand{\qbinn}[2]{\genfrac{[}{]}{0pt}{}{#1}{#2}_{\qpar}}

\newcommand{\pivo}{\ast}
\newcommand{\flip}{\star}

\newcommand{\gln}[1][n]{\mathfrak{gl}_{#1}}
\newcommand{\sln}[1][n]{\mathfrak{sl}_{#1}}

\newcommand{\webq}[1][\gln]{\catstuff{W}_{\qpar}{#1}}
\newcommand{\uwebq}[1][\gln]{\catstuff{W}_{\qpar}^{\uparrow}{#1}}

\newcommand{\uob}[1][k]{\uparrow_{#1}}
\newcommand{\dob}[1][k]{\downarrow_{#1}}

\newcommand{\shift}{d}
\newcommand{\formula}[1][{u,w}]{\chi(#1)}
\newcommand{\lformula}[1][{u,w}]{\chi(#1)}

\usepackage{aliascnt}
\def\NewTheorem#1{%
\newaliascnt{#1}{equation}%
\newtheorem{#1}[#1]{#1}%
\aliascntresetthe{#1}%
\expandafter\def\csname #1autorefname\endcsname{#1}%
}
\def\equationautorefname~#1\null{(#1)\null}

\numberwithin{equation}{subsection}

\NewTheorem{Proposition}
\NewTheorem{Theorem}
\NewTheorem{Corollary}
\AtEndEnvironment{Corollary}{\null\hfill$\square$}%
\NewTheorem{Lemma}
\theoremstyle{definition}
\NewTheorem{Definition}
\NewTheorem{Notation}
\NewTheorem{Definition/Notation}
\NewTheorem{Example}
\AtEndEnvironment{Example}{\null\hfill$\Diamond$}%
\theoremstyle{remark}
\NewTheorem{Remark}
\NewTheorem{Assumption}
\NewTheorem{Hypothesis}
\NewTheorem{Conjecture}
\NewTheorem{Question}
\NewTheorem{Observation}
\NewTheorem{Fact}
\NewTheorem{Conclusion}
\NewTheorem{lemma}

\usepackage[hypertexnames=false]{hyperref}
\usepackage{bookmark}
\hypersetup{
pdftoolbar=true,
pdfmenubar=true,
pdffitwindow=false,
pdfstartview={FitH},
pdftitle={A formula to evaluate type A webs and link polynomials},
pdfauthor={Abel Lacabanne, Daniel Tubbenhauer and Pedro Vaz},
pdfsubject={},
pdfcreator={Abel Lacabanne, Daniel Tubbenhauer and Pedro Vaz},
pdfproducer={Abel Lacabanne, Daniel Tubbenhauer and Pedro Vaz},
pdfkeywords={},
pdfnewwindow=true,
colorlinks=true,
linkcolor=mydarkblue,
citecolor=teal,
filecolor=magenta,
urlcolor=orchid,
linkbordercolor=lava,
citebordercolor=teal,
urlbordercolor=orchid,
linktocpage=true
}

\newcommand{\nnfootnote}[1]{%
\begin{NoHyper}
\renewcommand\thefootnote{}\footnote{#1}%
\addtocounter{footnote}{-1}%
\end{NoHyper}
}

\def\makeautorefname#1#2{\csdef{#1autorefname}{#2}}

\makeautorefname{section}{Section}%
\makeautorefname{subsection}{Section}%
\makeautorefname{subsubsection}{Section}%

\begin{document}
\title[A formula to evaluate type A webs and link polynomials]
{A formula to evaluate type A webs and link polynomials}
\author[A. Lacabanne, D. Tubbenhauer and P. Vaz]{Abel Lacabanne, Daniel Tubbenhauer and Pedro Vaz}

\address{A.L.: Laboratoire de Math{\'e}matiques Blaise Pascal (UMR 6620), Universit{\'e} Clermont Auvergne, Campus Universitaire des C{\'e}zeaux, 3 place Vasarely, 63178 Aubi{\`e}re Cedex, France,\newline \href{http://www.normalesup.org/~lacabanne}{www.normalesup.org/$\sim$lacabanne},
\href{https://orcid.org/0000-0001-8691-3270}{ORCID 0000-0001-8691-3270}}
\email{abel.lacabanne@uca.fr}

\address{D.T.: The University of Sydney, School of Mathematics and Statistics F07, Office Carslaw 827, NSW 2006, Australia, \href{http://www.dtubbenhauer.com}{www.dtubbenhauer.com}, \href{https://orcid.org/0000-0001-7265-5047}{ORCID 0000-0001-7265-5047}}
\email{daniel.tubbenhauer@sydney.edu.au}

\address{P.V.: Institut de Recherche en Math{\'e}matique et Physique, 
Universit{\'e} catholique de Louvain, Chemin du Cyclotron 2,  
1348 Louvain-la-Neuve, Belgium, \href{https://perso.uclouvain.be/pedro.vaz}{https://perso.uclouvain.be/pedro.vaz}, \href{https://orcid.org/0000-0001-9422-4707}{ORCID 0000-0001-9422-4707}}
\email{pedro.vaz@uclouvain.be}

\begin{abstract}
We give a closed formula to evaluate exterior webs (also called MOY webs) and the associated Reshetikhin--Turaev link polynomials.
\end{abstract}

\nnfootnote{\textit{Mathematics Subject Classification 2020.} Primary: 18M15; Secondary: 17B37, 18N25, 57K14.}
\nnfootnote{\textit{Keywords.} Webs and spiders, link polynomials, categorical skew Howe duality, web algebras, KLR algebras.}

\maketitle

\tableofcontents

\arrayrulewidth=0.5mm
\setlength{\arrayrulewidth}{0.5mm}


\section{Introduction}\label{S:Introduction}


The Reshetikhin--Turaev invariants of links \cite{ReTu-invariants-3-manifolds-qgroups} provide one of the most 
important family of link invariants in quantum topology and its ramifications. 

In this note we focus on the 
subfamily of these invariants given by coloring the strands 
of links with exterior powers of the vector representation of 
quantum $\sln$ or quantum $\gln$. We show that this 
family of polynomials can be computed by a \emph{closed formula} 
that takes as input only combinatorial data associated to 
a fixed colored link diagram and root theoretic data associated
to the type A Dynkin diagram. The formula works for all links, all exterior colorings, and all ranks $n$.

Let us stress that, by its very nature, the closed formula we give is 
a Weyl-character-type formula: On the one hand, 
it is general, completely algorithmic 
and might reveal abstract properties of the family of exterior colored link 
polynomials. 
But on the other hand, the formula does not necessarily give an efficient way to compute these invariants, at least not without massaging it a bit.

The diagrammatic incarnation of this family of link invariants is 
given by (exterior $\sln$ or $\gln$) \emph{webs} as shown in many 
works, see {\eg} \cite{RuTeWe-sl2}, \cite{Ku-spiders-rank-2}, \cite{MuOhYa-webs} or \cite{CaKaMo-webs-skew-howe}.
In fact, the aforementioned formula is an application of 
a \emph{closed evaluation formula} for webs 
that we will also state and prove. 
To use this formula one does not need to know any combinatorics of webs or tableaux; one, in fact, does not even need to know webs.
Moreover, another application of this formula 
is an easy to check criterion for webs to represent dual canonical basis elements. 

The exposition in this note, including the statement of 
our main results, is mostly 
self-contained and explicit. For example, we included 
Python based code, {\cf} \autoref{R:IntroductionEasyToCompute} below, that can do computations using, for example, the online calculator of SageMath. The main proofs, that are not necessary to understand the rest of the paper, however use techniques from categorification 
as we elaborate on now.

In breakthrough work Hu--Shi \cite{HuSh-monomial-klr-basis} 
found a closed formula for the dimension of the cyclotomic 
KLR algebra of any symmetrizable Kac--Moody type by using combinatorics of Fock spaces. Even more remarkable, 
their formula is a Weyl-character-type formula that can be computed without any prior knowledge of cyclotomic KLR algebras.
A consequence of categorical skew Howe duality is that 
cyclotomic KLR algebras of type A and webs (web algebras to be precise -- skew Howe duality itself it not enough to make the connection) are essentially the same object, see {\eg} the pioneering works \cite{BrSt-arc-algebras-1}, \cite{BrSt-arc-algebras-2} and \cite{BrSt-arc-algebras-3} 
for the $\sln[2]$ version of this result. See also \cite{MaPaTu-sl3-web-algebra} and \cite{Tu-intermediate-bases} for the $\sln[3]$ version, 
and \cite{Ma-sln-web-algebras} 
and \cite{Tu-gln-bases} for the general version of this relationship.

In this paper we put both together and obtain the aforementioned closed formula 
for web evaluations and the computation of link polynomials.

\begin{Remark}\label{R:IntroductionEasyToCompute}
The main formula in \autoref{Eq:EvaluationMain} is easy to 
compute via a machine, and so is its adjustment to the case 
of the link polynomials from \autoref{Eq:LinksMain}.
The reader may find Python code for SageMath
that can do the calculations 
here \cite{MaTu-sagemath-finite-type-klrw}. (The second link in \cite{MaTu-sagemath-finite-type-klrw} has 
``klrdim-webevaluation.py'' which is the version used 
for this paper. The two main examples 
\autoref{E:EvaluationMain} and \autoref{E:LinksTheInvariant} 
are included in that file. This code is also commented at the end of the source file.)
\end{Remark}

\begin{Remark}\label{R:IntroductionProofs}
To not distract the reader's attention, we postpone all proofs to 
\autoref{S:Proofs}. This has the advantage that we can formulate 
the main formulas without any reference to categorical skew Howe duality 
or KLR algebras.
\end{Remark}
\medskip

\noindent\textbf{Acknowledgments.}
We thank Jieru Zhu for helpful and related discussions, and Andrew Mathas 
for crucial help with the SageMath code which is also part of a joint project with D.T. (If the reader cites this code please also mention Andrew Mathas.)
We also thank the referee for a meticulous reading of the paper, and very helpful comments and suggestions.

This note is part of a project of D.T. trying to generalize 
\cite{BrSt-arc-algebras-1}, \cite{BrSt-arc-algebras-2}, \cite{BrSt-arc-algebras-3}, \cite{BrSt-arc-algebras-4}, and they have received much support along the way. Most recently, they like to 
thank a bike of the MFO Oberwolfach for breaking while riding: during 
the enforced rest period the main idea 
underlying this work was discovered.

A.L. was supported by a PEPS JCJC grant from INSMI (CNRS).
D.T. was supported by the Australian Research Council.
P.V. was supported by the Fonds de la Recherche Scientifique - FNRS under Grant no. MIS-F.4536.19.
We are also grateful to Universit{\'e} catholique de Louvain for supporting a research visit of D.T. to Louvain-la-Neuve.


\section{A reminder on webs}\label{S:Webs}


Fix $n\in\Z_{\geq 1}$.
We will now recall the description of (exterior $\gln$) webs 
from \cite{CaKaMo-webs-skew-howe}
using the exposition from \cite{LaTuVa-annular-webs-levi} and \cite{LaTu-gln-webs}. Background on monoidal categories 
given by generators-relations can be found in {\eg} 
\cite{TuVi-monoidal-tqft}, and background on monoidal categories related to our setting can be found in {\eg} \cite{EtGeNiOs-tensor-categories}.

We start with our reading and other conventions:

\begin{Notation}\label{N:WebsNotation}
All categories that we use are strict, and we read diagrams from bottom to top and left to right. The illustration
\begin{gather*}
(\idmor{}\hcirc\morstuff{g})\vcirc(\morstuff{f}\hcirc\idmor{}) 
=
\begin{tikzpicture}[anchorbase,scale=0.21,tinynodes]
\draw[spinach!35,fill=spinach!35] (-5.5,2) rectangle (-2,0.5);
\draw[spinach!35,fill=spinach!35] (-0.5,3) rectangle (3,4.5);
\draw[thick,densely dotted] (-5.5,2.5) node[left,yshift=-2pt]{$\vcirc$} 
to (3.5,2.5) node[right,yshift=-2pt]{$\vcirc$};
\draw[thick,densely dotted] (-1.25,5) 
node[above,yshift=-2pt]{$\hcirc$} to (-1.25,0) node[below]{$\hcirc$};
\draw[usual] (-5,2) to (-5,3.75) node[right]{$\dots$} to (-5,5);
\draw[usual] (-2.5,2) to (-2.5,5);
\draw[usual] (0,0) to (0,1.25) node[right]{$\dots$} to (0,3);
\draw[usual] (2.5,0) to (2.5,3);
\draw[usual] (-5,0) to (-5,0.25) node[right]{$\dots$} to (-5,0.5);
\draw[usual] (-2.5,0) to (-2.5,0.5);
\draw[usual] (2.5,4.5) to (2.5,5);
\draw[usual] (0,4.5) to (0,4.75) node[right]{$\dots$} to (0,5);
\node at (-3.75,1.0) {$\morstuff{f}$};
\node at (1.25,3.5) {$\morstuff{g}$};
\end{tikzpicture}
=
\begin{tikzpicture}[anchorbase,scale=0.21,tinynodes]
\draw[spinach!35,fill=spinach!35] (-5.5,1.75) rectangle (-2,3.25);
\draw[spinach!35,fill=spinach!35] (-0.5,1.75) rectangle (3,3.25);
\draw[usual] (-5,3.25) to (-5,4.125) node[right]{$\dots$} to (-5,5);
\draw[usual] (-2.5,3.25) to (-2.5,5);
\draw[usual] (0,0) to (0,0.625) node[right]{$\dots$} to (0,1.75);
\draw[usual] (2.5,0) to (2.5,1.75);
\draw[usual] (-5,0) to (-5,0.625) node[right]{$\dots$} to (-5,1.75);
\draw[usual] (-2.5,0) to (-2.5,1.75);
\draw[usual] (2.5,3.25) to (2.5,5);
\draw[usual] (0,3.25) to (0,4.125) node[right]{$\dots$} to (0,5);
\node at (-3.75,2.25) {$\morstuff{f}$};
\node at (1.25,2.25) {$\morstuff{g}$};
\end{tikzpicture}
=
\begin{tikzpicture}[anchorbase,scale=0.21,tinynodes]
\draw[spinach!35,fill=spinach!35] (5.5,2) rectangle (2,0.5);
\draw[spinach!35,fill=spinach!35] (0.5,3) rectangle (-3,4.5);
\draw[thick, densely dotted] (5.5,2.5) 
node[right,yshift=-2pt]{$\vcirc$} to (-3.5,2.5) node[left,yshift=-2pt]{$\vcirc$};
\draw[thick, densely dotted] (1.25,5) 
node[above,yshift=-2pt]{$\hcirc$} to (1.25,0) node[below]{$\hcirc$};
\draw[usual] (2.5,2) to (2.5,3.75) node[right]{$\dots$} to (2.5,5);
\draw[usual] (5,2) to (5,5);
\draw[usual] (0,0) to (0,3);
\draw[usual] (-2.5,0) to (-2.5,1.25) node[right]{$\dots$} to (-2.5,3);
\draw[usual] (2.5,0) to (2.5,0.25) node[right]{$\dots$} to (2.5,0.5);
\draw[usual] (5,0) to (5,0.5);
\draw[usual] (-2.5,4.5) to (-2.5,4.75) node[right]{$\dots$} to (-2.5,5);
\draw[usual] (0,4.5) to (0,5);
\node at (3.75,1.0) {$\morstuff{g}$};
\node at (-1.25,3.5) {$\morstuff{f}$};
\end{tikzpicture}
=(\morstuff{f}\hcirc\idmor{})\vcirc(\idmor{}\hcirc\morstuff{g}),
\end{gather*}
summarizes our reading conventions. The identity on an object $X$ is denoted by $\idmor{X}$.

The webs we will use below are certain labeled and oriented graphs.
The labels and orientations arise in 
a predetermined way from few choices. Using this we tend to omit 
many of the labels and orientations.

Although the edges of webs are labeled with $a\in\N$, we will 
allow labels $a\in\Z$ in formulas: by convention, the webs with negative labels are zero altogether.

We write $\munit$ for the monoidal unit
and $(\placeholder)^{\pivo}$ for the duality in 
a pivotal category. For webs $\munit$ is the empty word and 
$(\placeholder)^{\pivo}$ is changing upward to downward orientations for objects, extended to monoidal products using $(X\hcirc Y)^{\pivo}=Y^{\pivo}\hcirc X^{\pivo}$, 
and turning pictures by 180 degrees for morphisms.
\end{Notation}

Let $\qpar$ be a generic parameter.
For $a\in\Z$ and $b\in\N$ we let $\qnum{0}=0$, $\qfac{0}=1=\qbin{a}{0}$, $\qnum{a}=-\qnum{-a}$ for $a<0$ and 
otherwise
\begin{gather*}
\scalebox{0.97}{$\qnum{a}=\qpar^{a-1}+\qpar^{a-3}+\dots+\qpar^{-a+3}+\qpar^{-a+1}
,\quad
\qfac{b}=\qnum{b}\qnum{b-1}\dots\qnum{1}
,\quad
{\displaystyle\qbinn{a}{b}=
\frac{\qnum{a}\qnum{a-1}\dots\qnum{a-b+1}}{\qfac{b}}}$}.
\end{gather*}
The web categories we study are:

\begin{Definition}\label{D:WebsGln}
The \emph{(exterior $\gln$) web category $\webq$} is the
pivotal $\Zq$-linear category with $\hcirc$-generating objects 
$\uob$ and $\dob$, for $k\in\N$, 
of categorical dimension $\qbin{n}{k}$
and $\dob=(\uob)^{\pivo}$.

We further assume that $\webq$ has a braid group action on upwards objects, meaning morphisms depicted
\begin{gather*}
\text{over}:
\begin{tikzpicture}[anchorbase,scale=1]
\draw[usual,directed=0.99] (1,0)node[below]{$l$} to (0,1)node[above]{$l$};
\draw[usual,directed=0.99,crossline] (0,0)node[below]{$k$} to (1,1)node[above]{$k$};
\end{tikzpicture}
\colon\uob\hcirc\uob[l]\to\uob[l]\hcirc\uob
,\quad
\text{under}:
\begin{tikzpicture}[anchorbase,scale=1]
\draw[usual,directed=0.99] (0,0)node[below]{$k$} to (1,1)node[above]{$k$};
\draw[usual,directed=0.99,crossline] (1,0)node[below]{$l$} to (0,1)node[above]{$l$};
\end{tikzpicture}
\colon\uob\hcirc\uob[l]\to\uob[l]\hcirc\uob
\end{gather*}
for each simple braid group generator that 
satisfy the braid relations. We call these 
\emph{$(k,l)$-crossings} (overcrossings and undercrossings).
In these and similar pictures we tend to place the label under or over edges.

Finally, the $\vcirc$-$\hcirc$-generating morphisms are the ones coming from the pivotal structure and
\begin{gather*}
\scalebox{0.95}{$\begin{tikzpicture}[anchorbase,scale=1,yscale=-1]
\draw[usual,directed=0.99] (0.5,0.5) to (0,0)node[above]{$k$};
\draw[usual,directed=0.99] (0.5,0.5) to (1,0)node[above]{$l$};
\draw[usual,directed=0.99] (0.5,1)node[below]{$k{+}l$} to (0.5,0.5);
\end{tikzpicture}
\hspace*{-0.1cm}\colon
\uob[k{+}l]\to\uob\hcirc\uob[l]
,
\begin{tikzpicture}[anchorbase,scale=1,yscale=-1]
\draw[usual,directed=0.9] (0,0)node[above]{$k$} to (0.5,0.5);
\draw[usual,directed=0.9] (1,0)node[above]{$l$} to (0.5,0.5);
\draw[usual,directed=0.99] (0.5,0.5) to (0.5,1)node[below]{$k{+}l$};
\end{tikzpicture}
\hspace*{-0.1cm}\colon
\dob[k{+}l]\to\dob\hcirc\dob[l]
,
\begin{tikzpicture}[anchorbase,scale=1]
\draw[usual,directed=0.9] (0,0)node[below]{$k$} to (0.5,0.5);
\draw[usual,directed=0.9] (1,0)node[below]{$l$} to (0.5,0.5);
\draw[usual,directed=0.99] (0.5,0.5) to (0.5,1)node[above]{$k{+}l$};
\end{tikzpicture}
\hspace*{-0.1cm}\colon
\uob\hcirc\uob[l]\to\uob[k{+}l]
,
\begin{tikzpicture}[anchorbase,scale=1]
\draw[usual,directed=0.99] (0.5,0.5) to (0,0)node[below]{$k$};
\draw[usual,directed=0.99] (0.5,0.5) to (1,0)node[below]{$l$};
\draw[usual,directed=0.99] (0.5,1)node[above]{$k{+}l$} to (0.5,0.5);
\end{tikzpicture}
\hspace*{-0.1cm}\colon
\dob\hcirc\dob[l]\to\dob[k{+}l]
.$}
\end{gather*}

The relations imposed on $\webq$ are \emph{isotopies},
the \emph{exterior relation}, \emph{associativity}, \emph{coassociativity}, \emph{digon removal}, and 
\emph{dumbbell-crossing relation} together with \emph{invertibility} 
of the left mates of the $(k,l)$-overcrossings:
we take the quotient by the $\vcirc$-$\hcirc$-ideal generated by isotopies, invertibility and
\begin{gather*}
\begin{tikzpicture}[anchorbase,scale=1]
\draw[usual] (0,0)node[below]{${>}n$} to (0,1.5)node[above]{${>}n$};
\end{tikzpicture}
\hspace*{-0.15cm}=0
,
\begin{tikzpicture}[anchorbase,scale=1]
\draw[usual] (0,0)node[below]{$k$} to (1,1);
\draw[usual] (1,0)node[below]{$l$} to (0.5,0.5);
\draw[usual] (2,0)node[below]{$m$} to (1,1);
\draw[usual,directed=0.99] (1,1) to (1,1.5)node[above]{$k{+}l{+}m$};
\end{tikzpicture}
\hspace*{-0.3cm}=\hspace*{-0.3cm}
\begin{tikzpicture}[anchorbase,scale=1]
\draw[usual] (0,0)node[below]{$k$} to (1,1);
\draw[usual] (1,0)node[below]{$l$} to (1.5,0.5);
\draw[usual] (2,0)node[below]{$m$} to (1,1);
\draw[usual,directed=0.99] (1,1) to (1,1.5)node[above]{$k{+}l{+}m$};
\end{tikzpicture}
,
\begin{tikzpicture}[anchorbase,scale=1]
\draw[usual,rdirected=0.05] (0,0)node[above]{$k$} to (1,-1);
\draw[usual,rdirected=0.1] (1,0)node[above]{$l$} to (0.5,-0.5);
\draw[usual,rdirected=0.05] (2,0)node[above]{$m$} to (1,-1);
\draw[usual] (1,-1) to (1,-1.5)node[below]{$k{+}l{+}m$};
\end{tikzpicture}
\hspace*{-0.3cm}=\hspace*{-0.3cm}
\begin{tikzpicture}[anchorbase,scale=1]
\draw[usual,rdirected=0.05] (0,0)node[above]{$k$} to (1,-1);
\draw[usual,rdirected=0.1] (1,0)node[above]{$l$} to (1.5,-0.5);
\draw[usual,rdirected=0.05] (2,0)node[above]{$m$} to (1,-1);
\draw[usual] (1,-1) to (1,-1.5)node[below]{$k{+}l{+}m$};
\end{tikzpicture}
,
\begin{tikzpicture}[anchorbase,scale=1,rounded corners]
\draw[usual] (0.5,0.35) to (0,0.75)node[left]{$k$} to (0.5,1.15);
\draw[usual] (0.5,0.35) to (1,0.75)node[right]{$l$} to (0.5,1.15);
\draw[usual,directed=0.99] (0.5,1.15) to (0.5,1.5)node[above]{$k{+}l$};
\draw[usual] (0.5,0.35) to (0.5,0)node[below]{$k{+}l$};
\end{tikzpicture}
=\qbinn{k+l}{k}
\cdot\hspace*{-0.15cm}
\begin{tikzpicture}[anchorbase,scale=1]
\draw[usual,directed=0.99] (0.5,0)node[below]{$k{+}l$} to (0.5,1.5)node[above]{$k{+}l$};
\end{tikzpicture}
,\\
\begin{tikzpicture}[anchorbase,scale=1]
\draw[usual,directed=0.99] (0,0)node[below]{$k$} to (0.5,0.5) to (0.5,1) to (0,1.5)node[above]{$r$};
\draw[usual,directed=0.99] (1,0)node[below]{$l$} to (0.5,0.5) to (0.5,1) to (1,1.5)node[above]{$s$};
\end{tikzpicture}
=
(-1)^{kl}
\sum_{k-r=a-b}
(-\qpar)^{(k-a)(l-b)}\cdot
\begin{tikzpicture}[anchorbase,scale=1]
\draw[usual] (1,0.4)node[left,yshift=-0.05cm]{$b$} to (0,1.1);
\draw[usual,crossline] (0,0.4)node[right,yshift=-0.1cm]{$a$} to (1,1.1);
\draw[usual,directed=0.99] (0,0)node[below]{$k$} to (0,1.5)node[above]{$r$};
\draw[usual,directed=0.99] (1,0)node[below]{$l$} to (1,1.5)node[above]{$s$};
\end{tikzpicture}
=
(-1)^{kl}\sum_{k-r=a-b}
(-\qpar)^{-(k-a)(l-b)}\cdot
\begin{tikzpicture}[anchorbase,scale=1]
\draw[usual] (0,0.4)node[right,yshift=-0.1cm]{$a$} to (1,1.1);
\draw[usual,crossline] (1,0.4)node[left,yshift=-0.05cm]{$b$} to (0,1.1);
\draw[usual,directed=0.99] (0,0)node[below]{$k$} to (0,1.5)node[above]{$r$};
\draw[usual,directed=0.99] (1,0)node[below]{$l$} to (1,1.5)node[above]{$s$};
\end{tikzpicture}
,
\end{gather*}
together with their horizontally mirrored duals.
\end{Definition}

\begin{Remark}\label{R:WebsCrossings}
The final relation in \autoref{D:WebsGln} is also known as 
the \emph{Schur relation}. 
The non-quantum version comes from 
translating webs to the setting in Green's landmark book 
on the Schur algebra \cite{Gr-poly-reps}, hence the name. In \cite[Section 5]{LaTu-gln-webs} 
it is shown that this relation implies the more well-known crossing formula that we will use in 
\autoref{Eq:LinksSkein} below. 
(Strictly speaking \cite[Section 5]{LaTu-gln-webs} deals with symmetric webs, 
but the proof given therein works, {\muta}, for exterior webs.)
\end{Remark}

\begin{Notation}\label{N:WebsPhantom}
\leavevmode	

\begin{enumerate}

\item We also write $k$ for $\uob[k]$ and $-k$ for $\dob[k]$, where $k\in\N$. In this notation a general object of $\webq$ is of the form 
$\vec{k}=(k_{1},\dots,k_{r})\in\Z^{r}$ for some $r\in\N$.

\item A \emph{web} is a $\vcirc$-$\hcirc$-composition of the generating 
morphisms, {\ie} not a $\Zq$-linear combination.

\item There is no harm in thinking of webs as topological 
objects, meaning as labeled oriented graphs embedded in two-space.
We will sometimes use this to simplify drawings.

\item The edges labeled $n$ play a special role and we will illustrate them as
\begin{gather*}
\begin{tikzpicture}[anchorbase,scale=1,yscale=-1]
\draw[phantom,directed=0.99] (0,0.5) to (0,0)node[left,yshift=-0.25cm]{$n$};
\end{tikzpicture}
=
\begin{tikzpicture}[anchorbase,scale=1,yscale=-1]
\draw[usual,directed=0.99] (0,0.5) to (0,0)node[left,yshift=-0.25cm]{$n$};
\end{tikzpicture}
,\quad
\begin{tikzpicture}[anchorbase,scale=1,yscale=1]
\draw[phantom,directed=0.99] (0,0.5) to (0,0)node[left,yshift=0.25cm]{$n$};
\end{tikzpicture}
=
\begin{tikzpicture}[anchorbase,scale=1,yscale=1]
\draw[usual,directed=0.99] (0,0.5) to (0,0)node[left,yshift=0.25cm]{$n$};
\end{tikzpicture}
.
\end{gather*}
We call them \emph{phantom edges}. They should be thought of as nonexisting.
\end{enumerate}
\end{Notation}

It follows from the defining relations that the phantom edge calculus is essentially trivial, {\ie}:

\begin{Lemma}\label{L:WebsPhantom}
We have
\begin{gather*}
\begin{tikzpicture}[anchorbase,scale=1]
\draw[phantom,directed=0.99] (0,0) to[out=90,in=180] (0.5,0.5) to[out=0,in=90] (1,0)node[right]{$n$} to[out=270,in=0] (0.5,-0.5) to[out=180,in=270] (0,0);
\end{tikzpicture}
=1,\quad
\begin{tikzpicture}[anchorbase,scale=1]
\draw[phantom,directed=0.99] (1,1)node[above]{$n$} to (1,0)node[below]{$n$};
\draw[phantom,directed=0.99] (0,0)node[below]{$n$} to (0,1)node[above]{$n$};
\end{tikzpicture}
=
\begin{tikzpicture}[anchorbase,scale=1]
\draw[phantom,directed=0.55] (1,1)node[above]{$n$} to[out=270,in=270] (0,1)node[above]{$n$};
\draw[phantom,directed=0.55] (0,0)node[below]{$n$} to[out=90,in=90] (1,0)node[below]{$n$};
\end{tikzpicture}
,\quad
\begin{tikzpicture}[anchorbase,scale=1,rounded corners]
\draw[usual] (0.5,0.35) to (0,0.75)node[left]{$k$} to (0.5,1.15);
\draw[usual] (0.5,0.35) to (1,0.75)node[right]{$n{-}k$} to (0.5,1.15);
\draw[phantom,directed=0.99] (0.5,1.15) to (0.5,1.5)node[above]{$n$};
\draw[phantom] (0.5,0.35) to (0.5,0)node[below]{$n$};
\end{tikzpicture}
=
\begin{tikzpicture}[anchorbase,scale=1]
\draw[phantom,directed=0.99] (-0.5,-0.75)node[below]{$n$} to (-0.5,0.75)node[above]{$n$};
\draw[usual,directed=0.99] (0.299,0.25) to (0.3,0.25);
\draw[usual] (0,0) to[out=90,in=180] (0.25,0.25) to[out=0,in=90] (0.5,0)node[right]{$k$} to[out=270,in=0] (0.25,-0.25) to[out=180,in=270] (0,0);
\end{tikzpicture}
,\quad
\begin{tikzpicture}[anchorbase,scale=1]
\draw[phantom] (0,0)node[below]{$n$} to (0,0.4);
\draw[phantom,directed=0.99] (0,1.1) to (0,1.5)node[above]{$n$};
\draw[usual] (1,0)node[below]{$k$} to (1,0.6);
\draw[usual,directed=0.99] (1,0.9) to (1,1.5)node[above]{$k$};
\draw[usual] (1,0.9) to (0,1.1) to (0,0.4) to (1,0.6);
\end{tikzpicture}
=
\begin{tikzpicture}[anchorbase,scale=1]
\draw[phantom,directed=0.99] (0,0)node[below]{$n$} to (0,1.5)node[above]{$n$};
\draw[usual,directed=0.99] (1,0)node[below]{$k$} to (1,1.5)node[above]{$k$};
\end{tikzpicture}
.
\end{gather*}
There are more relations of a similar flavor which we omit to illustrate.
\end{Lemma}

Motivated by \autoref{L:WebsPhantom}, we call webs consisting 
of only phantom edges \emph{trivial}.

\begin{Remark}\label{R:WebsSlGl}
We will stay with $\webq$, which are webs for $\gln$, in this paper. 
We however stress that all results are valid 
for the respective 
$\sln$ version, thus, including the
Temperley--Lieb calculus \cite{RuTeWe-sl2}, 
Kuperberg's $\sln[3]$ spiders \cite{Ku-spiders-rank-2} and 
Cautis--Kamnitzer--Morrison's $\sln[n]$ webs 
\cite{CaKaMo-webs-skew-howe}. The (well-known) 
translation between these two pictures is 
a systematic identification of the form $\uob[k]\cong\dob[{n-k}]$.
This, representation theoretically, corresponds to 
the fact that $\bigwedge^{k}\C^{n}$ is dual to 
$\bigwedge^{n-k}\C^{n}$ as an $\sln$-module, 
but not as a $\gln$-module.
\end{Remark}

It will turn out to be useful to only allow \emph{upwards (pointing) webs}:

\begin{Definition}\label{D:FFormUpwards}
Let $\uwebq\subset\webq$ be the full subcategory
monoidally generated by $\{\uob[k]|k\in\N\}$.
\end{Definition}

We call webs in $\uwebq$ \emph{upwards webs}. By definition, webs in 
$\uwebq$ have upwards pointing boundary only, but can, a priori, have downwards oriented edges away from the boundary. One can show that $\uwebq$ 
has morphism spaces spanned by webs with all edges pointing upwards, but we will not need this fact.


\section{F-forms of webs}\label{S:FForm}


Let $\alpha_{i}=(0,\dots,0,1,-1,0,\dots,0)$, 
the \emph{$i$th simple root}, with the one in the $i$th entry.

\begin{Definition}\label{D:FFormFDef}
In this definition we work in $\uwebq$ only.
For all $a,i\in\N$ we define the \emph{$a$th F-operator} $F^{(a)}_{i}$ to be the operator that takes
$\idmor{\vec{k}}$, with $\vec{k}\in\N^{l}$, and returns $\idmor{\vec{k}-a\alpha_{i}}$, 
and we define 
the \emph{$a$th E-operator} $E^{(a)}_{i}$ to be the operator that takes
$\idmor{\vec{k}}$, with $\vec{k}\in\N^{l}$, and returns $\idmor{\vec{k}+a\alpha_{i}}$ 
given by
\begin{gather*}
\idmor{\vec{k}-a\alpha_{i}}F^{(a)}_{i}\idmor{\vec{k}}
=
\begin{tikzpicture}[anchorbase,scale=1]
\draw[usual] (1,0)node[below]{$k_{i{+}1}$} to (1,0.66);
\draw[usual,directed=0.99] (0,0.33) to (0,1)node[above,xshift=-0.1cm]{$k_{i}{-}a$};
\draw[usual,directed=0.99] (0,0)node[below]{$k_{i}$} to (0,0.33) to (1,0.66) to (1,1)node[above,xshift=0.1cm]{$k_{i{+}1}{+}a$};
\end{tikzpicture}
,\quad
\idmor{\vec{k}+a\alpha_{i}}E^{(a)}_{i}\idmor{\vec{k}}
=
\begin{tikzpicture}[anchorbase,scale=1]
\draw[usual] (1,0)node[below]{$k_{i{+}1}$} to (1,0.33);
\draw[usual,directed=0.99] (0,0.66) to (0,1)node[above,xshift=-0.1cm]{$k_{i}{+}a$};
\draw[usual,directed=0.99] (0,0)node[below]{$k_{i}$} to (0,0.66) to (1,0.33) to (1,1)node[above,xshift=0.1cm]{$k_{i{+}1}{-}a$};
\end{tikzpicture}
,
\end{gather*} 
and by the identity outside of these pictures.
The associated webs are called \emph{ladder webs}.
\end{Definition}

We will simplify notation using {\eg} $\idmor{\vec{k}-a\alpha_{i}}F^{(a)}_{i}\idmor{\vec{k}}=F^{(a)}_{i}\idmor{\vec{k}}=\idmor{\vec{k}-a\alpha_{i}}F^{(a)}_{i}$.
Moreover, for all $k\in\{1,\dots,n\}$ we use ladders to define
\begin{gather*}
\begin{tikzpicture}[anchorbase,scale=1,xscale=-1]
\draw[phantom,directed=0.99] (1,0)node[below]{$n$} to (0,1)node[above]{$n$};
\draw[usual,directed=0.99] (0,0)node[below]{$k$} to (1,1)node[above]{$k$};
\end{tikzpicture}
=
\begin{tikzpicture}[anchorbase,scale=1,xscale=-1]
\draw[phantom] (1,0)node[below]{$n$} to (1,0.33);
\draw[phantom,directed=0.99] (0,0.66) to (0,1)node[above]{$n$};
\draw[usual,directed=0.99] (0,0)node[below]{$k$} to (0,0.66) to (1,0.33) to (1,1)node[above]{$k$};
\end{tikzpicture}
,\quad
\begin{tikzpicture}[anchorbase,scale=1]
\draw[phantom,directed=0.99] (1,0)node[below]{$n$} to (0,1)node[above]{$n$};
\draw[usual,directed=0.99] (0,0)node[below]{$k$} to (1,1)node[above]{$k$};
\end{tikzpicture}
=
\begin{tikzpicture}[anchorbase,scale=1]
\draw[phantom] (1,0)node[below]{$n$} to (1,0.33);
\draw[phantom,directed=0.99] (0,0.66) to (0,1)node[above]{$n$};
\draw[usual,directed=0.99] (0,0)node[below]{$k$} to (0,0.66) to (1,0.33) to (1,1)node[above]{$k$};
\end{tikzpicture}
,\quad
\left(\text{note that }
\begin{tikzpicture}[anchorbase,scale=1]
\draw[phantom,directed=0.99] (1,0)node[below]{$n$} to (0,1)node[above]{$n$};
\draw[phantom,directed=0.99] (0,0)node[below]{$n$} to (1,1)node[above]{$n$};
\end{tikzpicture}
=
\begin{tikzpicture}[anchorbase,scale=1]
\draw[phantom,directed=0.99] (1,0)node[below]{$n$} to (1,1)node[above]{$n$};
\draw[phantom,directed=0.99] (0,0)node[below]{$n$} to (0,1)node[above]{$n$};
\end{tikzpicture}
\right)
.
\end{gather*}
We call these webs and all their mates \emph{phantom crossings}.
The following lemma allows us to 
use phantom crossings essentially without cost. Note that these are 
not the $(n,k)$-crossings coming from the braid group action on $\webq$, and in the picture above the strands ``cross virtually''. See \autoref{Eq:ProofsPhantom} for the precise relation between the various crossings.

\begin{Lemma}\label{L:FFormPhantom}
The phantom crossings satisfy all colored Reidemeister moves.
\end{Lemma}

\begin{Definition}\label{D:FFormPhantomForget}
We call all operations of the form
\begin{gather*}
\begin{tikzpicture}[anchorbase,scale=1,yscale=-1]
\draw[usual,directed=0.99] (0.5,0.5) to (0,0)node[above]{$k$};
\draw[usual,directed=0.99] (0.5,0.5) to (1,0)node[above]{$n-k$};
\draw[phantom] (0.5,1)node[below]{$n$} to (0.5,0.5);
\end{tikzpicture}
\mapsto
\begin{tikzpicture}[anchorbase,scale=1,yscale=-1]
\draw[white] (0.5,0.98)node[below,black]{$\munit$} to (0.5,0.5);
\draw[usual,directed=0.99] (0.5,0.5) to[out=180,in=90] (0,0)node[above]{$k$};
\draw[usual] (0.5,0.5) to[out=0,in=90] (1,0)node[above]{$k$};
\end{tikzpicture}
.
\end{gather*}
and all of its mates and mirrors \emph{forgetting phantom edges}. These operations can be successively applied to webs and will relabel and reorient them as part of these operations.

Recall from \autoref{N:WebsPhantom} that we think of webs as labeled oriented graphs. We say two such graphs are \emph{equivalent up to forgetting 
phantom edges} if one can be obtained from the other by
forgetting any finite number of phantom edges (including the relabeling and reorientation).
\end{Definition}

One can show that equivalence up to forgetting 
phantom edges defines an equivalence relation on webs.
By definition, webs are compositions of the generators of $\webq$, and we therefore can define:

\begin{Definition}\label{D:FFormUpwardsForm}
An \emph{upwards-form} of a web $w$ is a web $U(w)$ in $\uwebq$ 
that is equivalent to $w$ up to forgetting 
phantom edges.
\end{Definition}

\begin{Example}\label{E:FFormUpwardsForm}
Let $n=3$.
The web $w$ 
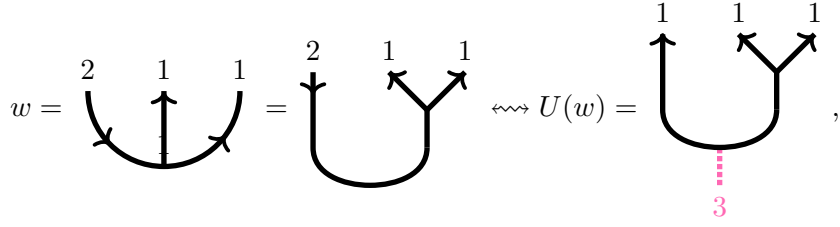
\begin{gather*}
w=
\begin{tikzpicture}[anchorbase,scale=1]
\draw[usual,directed=0.25,directed=0.8] (0,1)node[above]{$2$} to[out=270,in=180] (1,0) to[out=0,in=270] (2,1)node[above]{$1$};
\draw[usual,directed=0.99] (1,0) to (1,1)node[above]{$1$};
\end{tikzpicture}
=
\begin{tikzpicture}[anchorbase,scale=1]
\draw[usual,directed=0.1] (0,1)node[above]{$2$} to (0,0) to[out=270,in=180] (0.75,-0.5) to[out=0,in=270] (1.5,0);
\draw[usual] (1.5,0) to (1.5,0.5);
\draw[usual,directed=0.99] (1.5,0.5) to (2,1)node[above]{$1$};
\draw[usual,directed=0.99] (1.5,0.5) to (1,1)node[above]{$1$};
\end{tikzpicture}
\leftrightsquigarrow
U(w)=
\begin{tikzpicture}[anchorbase,scale=0.85]
\draw[phantom] (0.75,-1)node[below]{$3$} to (0.75,-0.5);
\draw[usual,rdirected=0.05] (0,1)node[above]{$1$} to (0,0) to[out=270,in=180] (0.75,-0.5) to[out=0,in=270] (1.5,0);
\draw[usual] (1.5,0) to (1.5,0.5);
\draw[usual,directed=0.99] (1.5,0.5) to (2,1)node[above]{$1$};
\draw[usual,directed=0.99] (1.5,0.5) to (1,1)node[above]{$1$};
\end{tikzpicture}
,
\end{gather*}
has the illustrated upwards-form.
\end{Example}

\begin{Lemma}\label{L:FFormUpwardsDef}
Every web has at least one upwards-form.
\end{Lemma}

\begin{Definition}\label{D:FFormDef}
An \emph{F-form} $F(w)=F_{i_{k}}^{(a_{r})}\dots F_{i_{1}}^{(a_{1})}\idmor{\vec{k}}$ of an upwards-pointing 
web $w$ is a string 
of $F$-operators such that the graphs of $w$ and
$F(w)$ are the same as labeled oriented graphs.

In general, an \emph{F-form} of a web $w$ is an F-form for any upwards-form $U(w)$.
\end{Definition}

\autoref{D:FFormDef} is best understood by example:

\begin{Example}\label{E:FFormDef}
For $n=3$, an F-form for $w$ as below is:
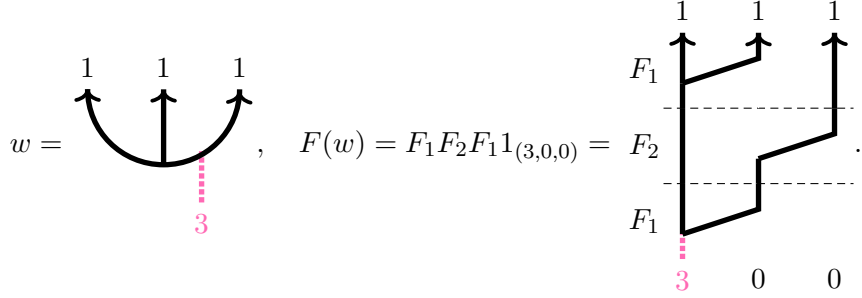
\begin{gather*}
w=
\begin{tikzpicture}[anchorbase,scale=1]
\draw[phantom] (1.5,-0.5)node[below]{$3$} to (1.5,0.18);
\draw[usual,directed=0.99] (0,0.99) to (0,1);
\draw[usual,directed=0.99] (0,1)node[above]{$1$} to[out=270,in=180] (1,0) to[out=0,in=270] (2,1)node[above]{$1$};
\draw[usual,directed=0.99] (1,0)
to (1,1)node[above]{$1$};
\end{tikzpicture}	
,\quad
F(w)=F_{1}F_{2}F_{1}\idmor{(3,0,0)}
=
\begin{tikzpicture}[anchorbase,scale=1]
\draw[white] (1,0)node[below,black]{$0$} to (1,0.66);
\draw[white] (2,0)node[below,black]{$0$} to (0,0.33);
\draw[usual] (0,0.33) to (0,1);
\draw[phantom] (0,0)node[below]{$3$} to (0,0.33);
\draw[usual] (0,0.33)to (1,0.66) to (1,1);
\draw[densely dashed] (-0.25,1) to (2.25,1);
\draw[usual] (0,1) to (0,2);
\draw[usual] (1,1) to (1,1.33);
\draw[usual] (1,1.33)to (2,1.66) to (2,2);
\draw[densely dashed] (-0.25,2) to (2.25,2);
\draw[white] (1,2) to (1,2.66);
\draw[usual,directed=0.99] (2,2) to (2,3)node[above]{$1$};
\draw[usual,directed=0.99] (0,2.33) to (0,2.66) to (0,3)node[above]{$1$};
\draw[usual] (0,2) to (0,2.33);
\draw[usual,directed=0.99] (0,2.33)to (1,2.66) to (1,3)node[above]{$1$};
\node at (-0.5,0.5) {$F_{1}$};
\node at (-0.5,1.5) {$F_{2}$};
\node at (-0.5,2.5) {$F_{1}$};
\end{tikzpicture}
.
\end{gather*}
Here and throughout, the horizontal slices are a visual aid only.
Note that F-forms are not unique and $F_{3}F_{1}F_{2}F_{1}\idmor{(3,0,0,0)}$ would be another F-form of $w$.
\end{Example}

We call objects of $\uwebq$ of the form $(n,\dots,n,0,\dots,0)=(n^{\ell},0,\dots,0)$ a \emph{level}.

\begin{Lemma}\label{L:FFormDef}
Every web has at least one F-form in $\Hom_{\uwebq}(\Lambda,\vec{k})$ 
for some level $\Lambda$.
\end{Lemma}

\begin{Definition}\label{D:FFormResidueSequence}
Let $I=\{i^{(a)}|i\in\Z_{\geq 1},a\in\N\}$.
Fix an F-form $F(w)=F_{i_{r}}^{(a_{r})}\dots F_{i_{1}}^{(a_{1})}\idmor{\Lambda}$ of $w$.
The \emph{residue sequence} for a web $w$ and its 
F-form is the tuple 
$r_{w}=(i_{1}^{(a_{1})},\dots,i_{r}^{(a_{r})})\in I^{r}$.
\end{Definition}

Note the reversed reading conventions when going from $F(w)$ to its residue sequence. 


\section{Evaluation of webs}\label{S:Evaluation}


Let ${}^{\flip}$ be the operation on webs that flips them upside down and reverses orientations.
We now define a pairing on $\Hom_{\webq}(\munit,\vec{k})$, which we call the \emph{evaluation pairing}:

\begin{Definition}\label{D:EvaluationPairing}
Given two webs $u,w\in\Hom_{\webq}(\munit,\vec{k})$ 
we let $(u,w)\in\End_{\webq}(\munit)$ be the element 
given by $(u,w)=w^{\flip}\vcirc u$, 
and we then extend this $\Zq$-linearly to all of 
$\Hom_{\webq}(\munit,\vec{k})$.
\end{Definition}

Our main goal is to give a closed formula for $(u,w)$. To this end, we need some preparation.

\begin{Definition}\label{D:EvaluationNCoefficients}
We define the following.	
\begin{enumerate}

\item Assume that we have an F-form 
$F(u)=F_{i_{r}}^{(a_{r})}\dots F_{i_{1}}^{(a_{1})}\idmor{\Lambda}$
of level $\Lambda=(n^{\ell},0,\dots,0)$ 
that ends at $\vec{k}=(k_{1},\dots,k_{m})$.
For $r_{u}=(i_{1}^{(a_{1})},\dots,i_{r}^{(a_{r})})$ we use the 
\emph{exploded sequence}
\begin{gather*}
\overline{r}_{u}
=
(\underbrace{i_{1},\dots,i_{1}}_{a_{1}
\text{ times}},\dots,\underbrace{i_{r},\dots,i_{r}}_{a_{r}
\text{ times}})
\text{ of length }\overline{r}=a_{1}+\dots+a_{r}.
\end{gather*}
We also let 
$\qfac{r_{u}}=\qfac{a_{1}}\dots\qfac{a_{r}}$. We will use $\tfrac{\qpar^{\shift}}{\qfac{r_{u}}\qfac{r_{w}}}$ which we call \emph{scaling}, where $\shift=\shift(\vec{k})=-\frac{1}{2}\big(n(n-1)\ell-\sum_{i=1}^{m}k_{i}(k_{i}-1)\big)$ is a shift.

\item We let $S_{\overline{r}}=\Aut(\{1,\dots,\overline{r}\})$ denote the symmetric group 
whose unit we denote by $e$.
For $\overline{r}_{u},\overline{r}_{w}\in I^{\overline{r}}$ we let $S_{\overline{r}_{u}}^{\overline{r}_{w}}=\{\sigma\in S_{\overline{r}}|\sigma\acts \overline{r}_{u}=\overline{r}_{w}\}$ be the \emph{set of possible 
crossings} where $\sigma\acts\overline{r}_{u}$ 
is the permutation of the entries of $\overline{r}_{u}$ determined by $\sigma$.

\item Recall the simple roots $\alpha_{i}$ and let
$\langle\alpha_{i},\alpha_{j}\rangle=a_{ij}$ with
be the usual Cartan 
pairing 
(so $a_{ii}=2$, $a_{ij}=-1$ for $|i-j|=1$ and 
$a_{ij}=0$ else) which we will use as indicated. 
Let $J_{\sigma}^{<t}=\{1\leq j<t|\sigma(j)<\sigma(t)\}$.
We define the following number of that we think as \emph{counting weights from previous steps}:
\begin{gather*}
N^{\Lambda}(\sigma,\overline{r}_{w},t)=
\Big\langle
\Lambda-\sum_{j\in J_{\sigma}^{<t}}\alpha_{(\overline{r}_{w})_{j}}
,\alpha_{(\overline{r}_{w})_{t}}
\Big\rangle
.
\end{gather*}
We also write $N(\sigma,t)=N^{\Lambda}(\sigma,\overline{r}_{w},t)$ for short.
	
\item We write $X_{\sigma}=\prod_{t=1}^{\overline{r}}
\qnum{N^{\Lambda}(\sigma,\overline{r}_{u},t)}\qpar^{N^{\Lambda}(e,\overline{r}_{u},t)-1}$. We let
\begin{gather}\label{Eq:EvaluationMain}
\formula=\tfrac{\qpar^{\shift}}{\qfac{r_{u}}\qfac{r_{w}}}\cdot
\sum_{\sigma\in S_{\overline{r}_{u}}^{\overline{r}_{w}}}
X_{\sigma}.
\end{gather}
	
\end{enumerate}
\end{Definition}

We encourage the reader to compare \autoref{Eq:EvaluationMain} to 
\cite[Theorem 1.1]{HuSh-monomial-klr-basis}. The setting of 
\cite{HuSh-monomial-klr-basis} looks quite different from ours, but 
we will see in \autoref{S:Proofs} why this comparison makes sense.

\begin{Theorem}\label{T:EvaluationMain}
For webs $u,w\in\Hom_{\webq}(\munit,\vec{k})$ we have that $\formula$ is independent of all choices involved, and
\begin{gather*}
(u,w)=\formula
\in\Nq
.
\end{gather*}
\end{Theorem}

\begin{Example}\label{E:EvaluationMain}
Readers are encouraged to refer to the SageMath code provided in Remark~\ref{R:IntroductionEasyToCompute} for further insight.

\begin{enumerate}

\item Let $n=2$. Let 
$u=w\in\Hom_{\webq}\big(\munit,(-1,1)\big)$ be, so that $(u,u)=\qnum{2}$:
\begin{gather*}
u=
\begin{tikzpicture}[anchorbase,scale=1]
\draw[white] (0,0) to (0,1);
\draw[usual,directed=0.99] (0,1)node[above]{$1$} to[out=270,in=180] (0.5,0.5) to[out=0,in=270] (1,1)node[above]{$1$};
\end{tikzpicture}
,\quad
(u,u)=u^{\flip}\vcirc u
=
\begin{tikzpicture}[anchorbase,scale=1]
\draw[white] (0,0) to (0,2);
\draw[usual,directed=0.99] (0,1) to[out=270,in=180] (0.5,0.5) to[out=0,in=270] (1,1);
\draw[usual] (0,1) to[out=90,in=180] (0.5,1.5) to[out=0,in=90] (1,1);
\end{tikzpicture}
=\qbinn{2}{1}=\qnum{2}.
\end{gather*}

As illustrated below, an F-form of $u$ is $F(u)=F_{1}\idmor{(2,0)}$, and 
an F-form of $u^{\flip}\vcirc u$ is $F(u^{\flip}\vcirc u)=F_{1}F_{1}\idmor{(2,0)}$.
We will also use the \emph{trivial web} with F-form $F(\emptyset)=F_{1}^{(2)}\idmor{(2,0)}$.
\begin{gather*}
F(u)=
\begin{tikzpicture}[anchorbase,scale=1]
\draw[white] (1,0)node[below,black]{$0$} to (1,0.66);
\draw[usual,directed=0.99] (0,0.33) to (0,1)node[above]{$1$};
\draw[phantom] (0,0)node[below]{$2$} to (0,0.33);
\draw[usual,directed=0.99] (0,0.33)to (1,0.66) to (1,1)node[above]{$1$};
\draw[densely dashed] (-0.25,1) to (1.25,1);
\node at (-0.5,0.5) {$F_{1}$};
\end{tikzpicture}
,\quad
F(u^{\flip}\vcirc u)=
\begin{tikzpicture}[anchorbase,scale=1]
\draw[white] (1,0)node[below,black]{$0$} to (1,0.66);
\draw[usual] (0,0.33) to (0,1);
\draw[phantom] (0,0)node[below]{$2$} to (0,0.33);
\draw[usual] (0,0.33)to (1,0.66) to (1,1);
\draw[densely dashed] (-0.25,1) to (1.25,1);
\draw[white] (0,2)node[above,black]{$0$} to (0,1.33);
\draw[phantom,directed=0.99] (1,1.66) to (1,2)node[above]{$2$};
\draw[usual] (1,1) to (1,1.66);
\draw[usual] (0,1) to (0,1.33);
\draw[usual] (0,1.33)to (1,1.66);
\node at (-0.5,0.5) {$F_{1}$};
\node at (-0.5,1.5) {$F_{1}$};
\end{tikzpicture}
,\quad
F(\emptyset)=
\begin{tikzpicture}[anchorbase,scale=1]
\draw[white] (1,0)node[below,black]{$0$} to (1,0.66);
\draw[white] (0,0.33) to (0,1)node[above,black]{$0$};
\draw[phantom,directed=0.99] (0,0)node[below]{$2$} to (0,0.33) to (1,0.66) to (1,1)node[above]{$2$};
\node at (-0.5,0.5) {$F_{1}^{(2)}$};
\end{tikzpicture}
.
\end{gather*}
The residue sequences from left to right are $(1)$, $(1,1)$ and $(1^{(2)})$, and in all cases $\Lambda=(2,0)$, 
and $\ell=1$, while $\vec{k}=(1,1)$ for $F(u)$ and $\vec{k}=(0,2)$ 
for the other two cases.

The main formula can be applied to either the pair $(u,u)$ or to 
$(\emptyset,u^{\flip}\vcirc u)$, giving the same result. That is, 
for $(u,u)$ we have $\overline{r}=1$,
$S_{\overline{r}_{u}}^{\overline{r}_{w}}=\{e\}$ and
\begin{gather*}
\sigma=e\colon X_{\sigma}=\underbrace{\qnum{2}\qpar}_{\qnum{N(e,1)}\qpar^{N(e,1)-1}}=1+\qpar^{2}.
\end{gather*}
Moreover, $\shift=-1$, no explosion was needed and
scaling by $\qpar^{-1}$ gives $\formula[{u,u}]=\qnum{2}$. 

For 
$(\emptyset,u^{\flip}\vcirc u)$ we have $\overline{r}=2$, $S_{\overline{r}_{u}}^{\overline{r}_{w}}=\{e,(1,2)\}$ and
\begin{gather*}
\sigma=e\colon X_{\sigma}=0
,\quad
\sigma=(1,2)\colon X_{\sigma}=
\hspace*{-1cm}
\underbrace{\qnum{2}\qpar}_{\qnum{N((1,2),1)}\qpar^{N(e,1)-1}}
\hspace*{-2cm}
\overbrace{\qnum{2}\qpar^{-1}}^{\qnum{N((1,2),2)}\qpar^{N(e,2)-1}}
\hspace*{-1cm}
=
\qnum{2}^{2}
.
\end{gather*}
This time we need to scale by $\qnum{2}^{-1}$ since 
we exploded $(1^{(2)})$ to $(1,1)$ during this computation and $\shift=0$. The result is the same.

\item Let $n=3$. We consider the following two webs $u,w\in\Hom_{\webq}\big(\munit,(-1,1,-1,1)\big)$:
\begin{gather*}
u=
\begin{tikzpicture}[anchorbase,scale=1]
\draw[white] (0,0) to (0,1);
\draw[usual,directed=0.99] (0,1)node[above]{$1$} to[out=270,in=180] (0.5,0.5) to[out=0,in=270] (1,1)node[above]{$1$};
\draw[usual,directed=0.99] (2,1)node[above]{$1$} to[out=270,in=180] (2.5,0.5) to[out=0,in=270] (3,1)node[above]{$1$};
\end{tikzpicture}
,\quad
w=
\begin{tikzpicture}[anchorbase,scale=1]
\draw[white] (0,0) to (0,1);
\draw[usual,directed=0.99] (0,1)node[above]{$1$} to[out=270,in=180] (1.5,0) to[out=0,in=270] (3,1)node[above]{$1$};
\draw[usual,directed=0.99] (2,1)node[above]{$1$} to[out=270,in=0] (1.5,0.5) to[out=180,in=270] (1,1)node[above]{$1$};
\end{tikzpicture}
.
\end{gather*}
We, of course, immediately get
\begin{gather*}
(u,w)=w^{\flip}\vcirc u
=
\begin{tikzpicture}[anchorbase,scale=1]
\draw[white] (0,0) to (0,2);
\draw[usual] (0,1) to[out=270,in=180] (0.5,0.5) to[out=0,in=270] (1,1);
\draw[usual,directed=0.99] (2,1) to[out=270,in=180] (2.5,0.5) to[out=0,in=270] (3,1);
\draw[usual] (0,1) to[out=90,in=180] (1.5,2) to[out=0,in=90] (3,1);
\draw[usual] (2,1) to[out=90,in=0] (1.5,1.5) to[out=180,in=90] (1,1);
\end{tikzpicture}
=\qbinn{3}{1}=\qnum{3}.
\end{gather*}
We get the same result from the main formula as follows. First, $F$-forms of $u$ and $w$ are
\begin{gather*}
F(u)=
\begin{tikzpicture}[anchorbase,scale=1]
\draw[white] (2,0)node[below,black]{$0$} to (2,0.66);
\draw[white] (3,0)node[below,black]{$0$} to (1,0.33);
\draw[usual] (1,0.33) to (1,1);
\draw[phantom] (0,0)node[below]{$3$} to (0,1);
\draw[phantom] (1,0)node[below]{$3$} to (1,0.33);
\draw[usual] (1,0.33)to (2,0.66) to (2,1);
\draw[densely dashed] (-0.25,1) to (3.25,1);
\draw[phantom] (0,1) to (0,2);
\draw[usual] (1,1) to (1,2);
\draw[usual] (2,1) to (2,1.33);
\draw[usual] (2,1.33)to (3,1.66) to (3,2);
\draw[densely dashed] (-0.25,2) to (3.25,2);
\draw[phantom] (0,2) to (0,3);
\draw[usual] (2,3) to (2,2.66) to (1,2.33) to (1,2);
\draw[usual] (3,2) to (3,3);
\draw[densely dashed] (-0.25,3) to (3.25,3);
\draw[phantom] (0,3) to (0,3.33);
\draw[usual,directed=0.99] (1,3.66) to (1,4)node[above]{$1$};
\draw[usual,directed=0.99] (1,3.66) to (0,3.33) to (0,4)node[above]{$2$};
\draw[usual,directed=0.99] (2,3) to (2,4)node[above]{$2$};
\draw[usual,directed=0.99] (3,3) to (3,4)node[above]{$1$};
\node at (-0.5,0.5) {$F_{2}$};
\node at (-0.5,1.5) {$F_{3}$};
\node at (-0.5,2.5) {$F_{2}^{(2)}$};
\node at (-0.5,3.5) {$F_{1}$};
\end{tikzpicture}
,\quad
F(w)=
\begin{tikzpicture}[anchorbase,scale=1]
\draw[white] (2,0)node[below,black]{$0$} to (2,0.66);
\draw[white] (3,0)node[below,black]{$0$} to (1,0.33);
\draw[usual] (1,0.33) to (1,1);
\draw[phantom] (0,0)node[below]{$3$} to (0,1);
\draw[phantom] (1,0)node[below]{$3$} to (1,0.33);
\draw[usual] (1,0.33)to (2,0.66) to (2,1);
\draw[densely dashed] (-0.25,1) to (3.25,1);
\draw[phantom] (0,1) to (0,2);
\draw[usual] (1,1) to (1,2);
\draw[usual] (2,1) to (2,1.33);
\draw[usual] (2,1.33)to (3,1.66) to (3,2);
\draw[densely dashed] (-0.25,2) to (3.25,2);
\draw[phantom] (0,2) to (0,2.33);
\draw[phantom] (1,2.66) to (1,3);
\draw[usual] (1,2) to (1,2.66) to (0,2.33) to (0,3);
\draw[usual] (3,2) to (3,3);
\draw[densely dashed] (-0.25,3) to (3.25,3);
\draw[usual,directed=0.99] (0,3) to (0,4)node[above]{$2$};
\draw[phantom] (1,3) to (1,3.33);
\draw[usual,directed=0.99] (2,3.66) to (2,4)node[above]{$2$};
\draw[usual,directed=0.99] (2,3.66) to (1,3.33) to (1,4)node[above]{$1$};
\draw[usual,directed=0.99] (3,3) to (3,4)node[above]{$1$};
\node at (-0.5,0.5) {$F_{2}$};
\node at (-0.5,1.5) {$F_{3}$};
\node at (-0.5,2.5) {$F_{1}$};
\node at (-0.5,3.5) {$F_{2}^{(2)}$};
\end{tikzpicture}
.
\end{gather*}
Here $\Lambda=(3^{2},0,0)$ and $\ell=2$.
Using the associated residue sequences $(2,3,2^{(2)},1)$ and 
$(2,3,1,2^{(2)})$ we get the same result: Firstly, the exploded sequences 
are $(2,3,2,2,1)$ and 
$(2,3,1,2,2)$ so we remember that we have to multiply in the end by $\qpar^{-4}\qnum{2}^{-2}$. We compute 
\begin{gather*}
S_{\overline{r}_{u}}^{\overline{r}_{w}}=
\{
(3,4,5), (1,4)(3,5), (1,5,3,4), (3,5), (1,4,5,3), (1,5,3)
\}
\end{gather*}
where we use the usual notation for permutations. The six relevant summands are then
\begin{gather*}
\sigma\neq(1,4,5,3)\colon X_{\sigma}=0
,\quad
\sigma=(1,4,5,3)\colon X_{\sigma}=
\qpar^{4}\qnum{2}^{2}\qnum{3}
=(1+\qpar^{2})^{2}(1+\qpar^{2}+\qpar^{4})
.
\end{gather*}
Thus, scaling by $\qpar^{-4}\qnum{2}^{-2}$ gives the desired result. Here 
we explode twice, so we get $\qnum{2}^{-2}$ and $\shift=-4$.

As in (a), we could also use an F-form for $w^{\flip}\vcirc u$, 
{\eg} $(2,3,2^{(2)},1,2,1^{(2)},3^{(2)},2^{(2)})$, which we pair with the trivial web $(2^{(3)},3^{(3)},1^{(3)},2^{(3)})$. Applying the formula for the exploded 
residue sequences $(2,3,2,2,1,2,1,1,3,3,2,2)$ and 
$(2,2,2,3,3,3,1,1,1,2,2,2)$ gives the same result.

\end{enumerate}
For completeness, let $(n^{\ell},0,\dots,0)$ have $\ell$ symbols $n$ and 
$m-\ell$ symbols $0$.
In general a residue sequence of the trivial web 
see as a web in $\Hom_{\webq}\big((n^{\ell},0,\dots,0),(0,\dots,0,n^{\ell})\big)$ is
given by the residue sequence $(\ell^{(n)},\dots,(m-1)^{(n)},(\ell-1)^{(n)},\dots,(m-2)^{(n)},\dots)$.
\end{Example}

Due to their relation to invariant tensors, the space 
$\Hom_{\webq}(\munit,\vec{k})$ has an important basis 
known as \emph{Lusztig--Kashiwara's dual canonical basis}. 
(For details see {\eg} \cite{KhKu-sl3-web-bases} or \cite[Section 4.1.5]{Tu-gln-bases}.)
We can thus ask whether a fixed $w\in\Hom_{\webq}(\munit,\vec{k})$
corresponds to a dual canonical basis element, and 
\autoref{T:EvaluationMain} gives a complete answer 
that can be checked using 
the main formula:

\begin{Proposition}\label{P:EvaluationDualCanonical}
A $w\in\Hom_{\webq}(\munit,\vec{k})$ is dual canonical if and only if 
$\qpar^{-\shift}\formula[{w,w}]\in\Nqq$.
\end{Proposition}


\section{Evaluation of link polynomials}\label{S:Links}


In order to compute link polynomials using F-forms we first need 
to explain how to interpret crossings. 
To this end, note that the category $\webq$ has no generating sideways 
or downwards crossings as extra generators. However, 
these crossings can be obtained as compositions 
of upwards-pointing crossings and the generators of $\webq$.
For example:
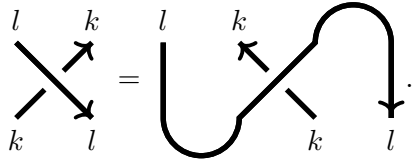
\begin{gather*}
\begin{tikzpicture}[anchorbase,scale=1]
\draw[usual,directed=0.99] (0,0)node[below]{$k$} to (1,1)node[above]{$k$};
\draw[usual,directed=0.99,crossline] (0,1)node[above]{$l$} to (1,0)node[below]{$l$};
\end{tikzpicture}
=
\begin{tikzpicture}[anchorbase,scale=1]
\draw[usual,directed=0.99] (1,0)node[below]{$k$} to (0,1)node[above]{$k$};
\draw[usual,directed=0.99,crossline] (-1,1)node[above]{$l$} to (-1,0) to[out=270,in=180] (-0.5,-0.5) to[out=0,in=270] (0,0) to (1,1) to[out=90,in=180] (1.5,1.5) to[out=0,in=90] (2,1) to (2,0)node[below]{$l$};
\end{tikzpicture}
.
\end{gather*}

We consider framed oriented colored links $L_{\vec{c}}$ in $\R^{3}$, whose components are colored by $\vec{c}\in\N^{r}$ where $r\in\N$ is the number of components.
We associate F-forms to their diagrams $D_{\vec{c}}$:

\begin{Definition}\label{D:LinksFLink}
Given an oriented link diagram $D_{\vec{c}}$. An \emph{F-form} of $D_{\vec{c}}$ is the 
linear combinations of webs obtained by:
\begin{enumerate}

\item Choosing a Morse positioning (in terms of the $\vcirc$-$\hcirc$ 
generators, including upwards-pointing crossings) of $D_{\vec{c}}$;

\item Replace all $(k,l)$-crossings by 
\begin{gather}\label{Eq:LinksSkein}
\begin{tikzpicture}[anchorbase,scale=1]
\draw[usual,directed=0.99] (1,0)node[below]{$l$} to (0,1)node[above]{$l$};
\draw[usual,directed=0.99,crossline] (0,0)node[below]{$k$} to (1,1)node[above]{$k$};
\end{tikzpicture}
=
(-1)^{kl}
\sum_{b-a=k-l}
(-\qpar)^{k-b}
\begin{tikzpicture}[anchorbase,scale=1]
\draw[usual,directed=0.99] (0,0)node[below]{$k$} to (0,1.5)node[above]{$l$};
\draw[usual,directed=0.99] (1,0)node[below]{$l$} to (1,1.5)node[above]{$k$};
\draw[usual] (0,0.4) to node[below]{$b$} (1,0.6);
\draw[usual] (1,0.9) to node[above]{$a$} (0,1.1);
\end{tikzpicture}
,\quad
\begin{tikzpicture}[anchorbase,scale=1]
\draw[usual,directed=0.99] (0,0)node[below]{$k$} to (1,1)node[above]{$k$};
\draw[usual,directed=0.99,crossline] (1,0)node[below]{$l$} to (0,1)node[above]{$l$};
\end{tikzpicture}
=
(-1)^{kl}
\sum_{b-a=k-l}
(-\qpar)^{-k+b}
\begin{tikzpicture}[anchorbase,scale=1]
\draw[usual,directed=0.99] (0,0)node[below]{$k$} to (0,1.5)node[above]{$l$};
\draw[usual,directed=0.99] (1,0)node[below]{$l$} to (1,1.5)node[above]{$k$};
\draw[usual] (0,0.4) to node[below]{$b$} (1,0.6);
\draw[usual] (1,0.9) to node[above]{$a$} (0,1.1);
\end{tikzpicture}
,
\end{gather}
and obtain a linear combinations of webs 
$a_{1}w_{1}+\dots+a_{s}w_{s}$ where $a_{i}\in\Zq$;

\item Replace all $w_{i}$ by a choice of F-form $F(w_{i})$;

\item The associated F-form of $D_{\vec{c}}$ is $F(D_{\vec{c}})=
a_{1}F(w_{1})+\dots+a_{s}F(w_{s})$.
	
\end{enumerate}
\end{Definition}

\begin{Lemma}\label{L:LinksFLink}
Every colored link diagram has at least one F-form.
\end{Lemma}

\begin{Definition}\label{L:LinksTheInvariant}
Let $L_{\vec{c}}$ be a framed oriented colored link 
and let $F(D_{\vec{c}})=
a_{1}F(w_{1})+\dots+a_{s}F(w_{s})$ be any F-form of any of its diagrams.
We define:
\begin{gather}\label{Eq:LinksMain}
\lformula[L_{\vec{c}}]
=\sum_{i=1}^{s}a_{i}\formula[{w_{i}}]
\in\Zq
.
\end{gather}
\end{Definition}

\begin{Remark}\label{R:LinksFramed}
For the reader who wants to work with links and not framed links: One can easily verify that a $k$-colored 
Reidemeister I move gives the scalars $\qpar^{\pm(k(-k+n+1))}$ 
(plus for the overcrossing and minus for the undercrossing). This in turn 
determines how the invariant in 
\autoref{L:LinksTheInvariant} changes under framing changes.
\end{Remark}

\begin{Example}\label{E:LinksTheInvariant}
This example can also be found in the SageMath code as in \autoref{R:IntroductionEasyToCompute}.

Let $n=2$ and let $H_{(1,1)}$ denote the $(1,1)$-colored 
Hopf link coming from the braid word with two positive crossings.
The reader may convince themselves that an F-form of the 
standard diagram $D_{(1,1)}$ of $H_{(1,1)}$ is
\begin{gather*}
\scalebox{0.93}{$F(D_{(1,1)})=
F_{4}F_{5}F_{3}F_{4}
(\qpar^{2}F_{2}F_{3}F_{1}F_{2}
-\qpar F_{2}F_{3}F_{2}F_{1}
-\qpar F_{3}F_{2}F_{1}F_{2}
+ F_{3}F_{2}F_{2}F_{1})
F_{4}F_{3}F_{5}F_{4}F_{2}F_{3}F_{1}F_{2}\idmor{(2,2,0,0,0,0)}$}
,
\end{gather*}
which we compare to the trivial web with residue sequence 
$(2^{(2)},3^{(2)},4^{(2)},5^{(2)},1^{(2)},2^{(2)},3^{(2)},4^{(2)})$. Using the 
main formula from \autoref{T:EvaluationMain} we get 
\begin{gather*}
(\emptyset,H_{(1,1)})=\qpar^{2}(\qpar+\qpar^{-1})^{2}-2\qpar(\qpar+\qpar^{-1})+(\qpar+\qpar^{-1})^{2}
=\qpar^{-1}\qnum{4},
\end{gather*}
which is the expected result up to a power of $\qpar$. (The precise power depends on the conventions 
one wants to compare the above to.)
\end{Example}

\begin{Theorem}\label{T:LinksTheInvariant}
The Laurent polynomial 
$\lformula[L_{\vec{c}}]$ is well-defined, {\ie} 
it is an invariant of $L_{\vec{c}}$ and independent of all choices involved.
Moreover, up to potentially adjusting conventions, this invariant agrees with the Reshetikhin--Turaev exterior colored $\gln$ link invariant.
\end{Theorem}

\begin{Remark}\label{R:LinksTheInvariant}
The reader who wants to work with $\sln$ instead of $\gln$ 
needs to shift the Laurent polynomial 
$\lformula[L_{\vec{c}}]$ in \autoref{T:LinksTheInvariant} by $\qpar^{1/n}$ to match {\eg} the Reshetikhin--Turaev exterior colored $\sln$ link invariant in \cite[Corollary 6.2.3]{CaKaMo-webs-skew-howe}.
\end{Remark}


\section{The proofs}\label{S:Proofs}


We now give all the proofs, sometimes collected into one proof.

\begin{proof}[Proof of ``The phantom calculus is trivial, {\ie} \autoref{L:WebsPhantom} and \autoref{L:FFormPhantom}'']
\textit{(Part a.)} All of the relations displayed in \autoref{L:WebsPhantom} 
follow directly from the definitions except the right-hand equation 
for which we use \cite[Lemma 2A.14]{LaTu-gln-webs}.

\textit{(Part b.)}  It follows from (the exterior version of) \cite[Section 5]{LaTu-gln-webs} that
\begin{gather}\label{Eq:ProofsPhantom}
\scalebox{0.95}{$\begin{tikzpicture}[anchorbase,scale=1,xscale=-1]
\draw[phantom,directed=0.99] (1,0)node[below]{$n$} to (0,1)node[above]{$n$};
\draw[usual,directed=0.99] (0,0)node[below]{$k$} to (1,1)node[above]{$k$};
\end{tikzpicture}
=
(-1)^{nk+k}\qpar^{-k}\cdot
\begin{tikzpicture}[anchorbase,scale=1,xscale=-1]
\draw[usual,directed=0.99] (0,0)node[below]{$k$} to (1,1)node[above]{$k$};
\draw[phantom,directed=0.99,crosslinep] (1,0)node[below]{$n$} to (0,1)node[above]{$n$};
\end{tikzpicture}
=
(-1)^{nk+k}\qpar^{k}\cdot
\begin{tikzpicture}[anchorbase,scale=1,xscale=-1]
\draw[phantom,directed=0.99] (1,0)node[below]{$n$} to (0,1)node[above]{$n$};
\draw[usual,directed=0.99,crossline] (0,0)node[below]{$k$} to (1,1)node[above]{$k$};
\end{tikzpicture}$}
,
\end{gather}
and similarly for the other phantom crossings.
(In the above picture note the difference between the honest $(n,k)$-crossings and the phantom crossings.) As a consequence of 
{\eg} (the exterior version of) 
\cite[Section 2]{LaTu-gln-webs}, the honest $(k,l)$-crossings 
satisfy the Reidemeister relations, with Reidemeister I only up to the scalar $\qpar^{\pm(k(-k+n+1))}$ (here $k=l$), and the other defining relations of the colored tangle 
category (see {\eg} for an uncolored list of these relations). Thus, 
the statement of \autoref{L:FFormPhantom} follows from \autoref{Eq:ProofsPhantom}.
\end{proof}

\begin{proof}[Proof of ``F-forms exist, {\ie} \autoref{L:FFormUpwardsDef}, \autoref{L:FFormDef} and \autoref{L:LinksFLink}'']
\textit{(Part a.)} That every web, more precisely any expression in the 
$\vcirc$-$\hcirc$ generators, has an upwards-form can be seen inductively: Let $h\in\N$ denote the number of Morse points. If $h=0$ and the web we start with is already upwards-pointing, then there is nothing to show. If the starting web is downwards pointing, then we can just reverse all orientations. The analog for $h=1$ is also easily verified.
So assume $h>1$. Pick any Morse point and perform 
either of
\begin{gather}\label{Eq:ProofPhantom}
\scalebox{0.85}{$\begin{tikzpicture}[anchorbase,scale=1,yscale=-1]
\draw[white] (0.5,0.98)node[below,black]{$\munit$} to (0.5,0.5);
\draw[usual,directed=0.99] (0.5,0.5) to[out=180,in=90] (0,0)node[above]{$k$};
\draw[usual] (0.5,0.5) to[out=0,in=90] (1,0)node[above]{$k$};
\end{tikzpicture}
\mapsto
\begin{tikzpicture}[anchorbase,scale=1,yscale=-1]
\draw[usual,directed=0.99] (0.5,0.5) to (0,0)node[above]{$k$};
\draw[usual,directed=0.99] (0.5,0.5) to (1,0)node[above]{$n-k$};
\draw[phantom] (0.5,1)node[below]{$n$} to (0.5,0.5);
\end{tikzpicture}
,\quad
\begin{tikzpicture}[anchorbase,xscale=-1,yscale=-1]
\draw[white] (0.5,0.98)node[below,black]{$\munit$} to (0.5,0.5);
\draw[usual,directed=0.99] (0.5,0.5) to[out=180,in=90] (0,0)node[above]{$k$};
\draw[usual] (0.5,0.5) to[out=0,in=90] (1,0)node[above]{$k$};
\end{tikzpicture}
\mapsto
\begin{tikzpicture}[anchorbase,xscale=-1,yscale=-1]
\draw[usual,directed=0.99] (0.5,0.5) to (0,0)node[above]{$k$};
\draw[usual,directed=0.99] (0.5,0.5) to (1,0)node[above]{$n-k$};
\draw[phantom] (0.5,1)node[below]{$n$} to (0.5,0.5);
\end{tikzpicture}
,\quad
\begin{tikzpicture}[anchorbase,xscale=1,yscale=1]
\draw[white] (0.5,0.98)node[above,black]{$\munit$} to (0.5,0.5);
\draw[usual,directed=0.99] (0.5,0.5) to[out=180,in=90] (0,0)node[below]{$k$};
\draw[usual] (0.5,0.5) to[out=0,in=90] (1,0)node[below]{$k$};
\end{tikzpicture}
\mapsto
\begin{tikzpicture}[anchorbase,xscale=-1,yscale=1]
\draw[usual] (0,0)node[below]{$k$} to (0.5,0.5);
\draw[usual] (1,0)node[below]{$n-k$} to (0.5,0.5);
\draw[phantom,directed=0.99] (0.5,0.5) to (0.5,1)node[above]{$n$};
\end{tikzpicture}
,\quad
\begin{tikzpicture}[anchorbase,xscale=-1,yscale=1]
\draw[white] (0.5,0.98)node[above,black]{$\munit$} to (0.5,0.5);
\draw[usual,directed=0.99] (0.5,0.5) to[out=180,in=90] (0,0)node[below]{$k$};
\draw[usual] (0.5,0.5) to[out=0,in=90] (1,0)node[below]{$k$};
\end{tikzpicture}
\mapsto
\begin{tikzpicture}[anchorbase,xscale=1,yscale=1]
\draw[usual] (0,0)node[below]{$k$} to (0.5,0.5);
\draw[usual] (1,0)node[below]{$n-k$} to (0.5,0.5);
\draw[phantom,directed=0.99] (0.5,0.5) to (0.5,1)node[above]{$n$};
\end{tikzpicture}$}
\end{gather}
once, if $h$ is odd, or twice (at two different Morse points), 
if $h$ is even.
Using phantom crossings we connected 
the dangling phantom strands in those pictures to anywhere 
point at the bottom and top, respectively. Relabeling and reorienting the 
result gives a legal web with smaller $h$ by the combinatorics 
of oriented plane trivalent graphs (we can ignore the precise positions of 
the phantom edges by \autoref{L:FFormPhantom}), so the claim of 
\autoref{L:FFormUpwardsDef} follows.

\textit{(Part b.)} This can be proven by using 
\autoref{L:FFormUpwardsDef} and \cite[Lemma 4.9]{Tu-gln-bases}.

\textit{(Part c.)} \autoref{L:LinksFLink} follows directly from \autoref{L:FFormDef} via the Skein relations \autoref{Eq:LinksSkein}.
\end{proof}

\begin{proof}[Proof of ``The formulas work, {\ie} \autoref{T:EvaluationMain} and \autoref{T:LinksTheInvariant}'']
\textit{(Part a.)} In the proof below we will use various statements about $\gln$ web algebras, all 
of which, more or less explicitly, can be found in 
\cite{Ma-sln-web-algebras} or \cite{Tu-gln-bases}, using a matrix factorization description.

For $u,w\in\Hom_{\webq}(\munit,\vec{k})$ 
let $W(u,w)$ denote the associated free $\Zq$-module summand 
of the $\gln$ web algebra, see \cite[Section 3.3.4]{Tu-gln-bases} for details. More precisely, $W(u,w)$ is the idempotent truncation of the 
$\gln$ web algebra obtained by using the two webs $u$ and $w$.
By the construction of the $\gln$ web algebra via $\gln$ foams or 
$\gln$ matrix factorizations, respectively, 
and the universal construction, 
we have that $(u,w)=\qpar^{\shift}
\mathrm{rank}_{\Zq}\big(W(u,w)\big)$ (the graded rank) and we will use 
this throughout. (This graded rank formula can be obtained in many ways. The paper \cite{RoWa-foam-formula} gives a self-contained summary 
how the universal construction applied to $\gln$ foams gives 
a categorification of the MOY calculus and thus the formula $(u,w)=\qpar^{\shift}
\mathrm{rank}_{\Zq}\big(W(u,w)\big)$ follows.) Here the shift by $\qpar^{\shift}$ comes simply 
from the desire to have the unit of the $\gln$ web algebra sitting in degree zero, which corresponds to using Gaussians ({\eg} $1+\qpar^{2}$) instead of quantum numbers.

Let us first assume that we have $u,w\in\Hom_{\uwebq}(\Lambda,\vec{k})$ for some level $\Lambda$ with fixed F-forms.
The result in \cite[Theorem 5.16]{Tu-gln-bases} identifies 
$W(u,w)$ as an idempotent truncation of a thick version of the 
cyclotomic KLR algebra of type $A_{\Z}$ and level $\Lambda$. 
The idempotent truncation is exactly given by the images 
of the F-forms of $u,w$ under the categorified skew Howe duality 
in terms of the $\gln$ web algebra, see {\eg} \cite[Section 3.3]{Tu-gln-bases}. The thickening, as explained in 
\cite[Section 3.3]{Tu-gln-bases} is just the KLR diagram 
version of the usual thick calculus from \cite{KhLaMaSt-thick-calculus}, 
and it is easy to see that explosion in this case corresponds to 
isomorphisms on both sides of categorified skew Howe duality. As 
explained in \cite[Section 5]{Tu-gln-bases} these various isomorphism 
patched together provide the factors $\qfac{r_{u}}$ and $\qfac{r_{w}}$ in the formula from \autoref{Eq:EvaluationMain}. Otherwise, \autoref{Eq:EvaluationMain} is the type $A_{\Z}$ version of
\cite[Theorem A]{HuSh-monomial-klr-basis}, so we can use \cite[Theorem 5.16]{Tu-gln-bases} to push it over to the $\gln$ web algebra.

It remains to argue that choosing different F-forms does not 
change the result. If $F(u)$ and $F(u)^{\prime}$ are two 
F-forms for $u$ that are equal as topological webs and similarly 
$F(w)$ and $F(w)^{\prime}$ are two 
F-forms for $u$ that are equal as topological webs, then 
$(F(u),F(w))=(F(u)^{\prime},F(w)^{\prime})$ using the topological 
invariance of the $\gln$ web algebra. By the previous point we 
further know that $(F(u),F(w))=\formula$ and $(F(u)^{\prime},F(w)^{\prime})=\formula[{u,w}]^{\prime}$, with the prime indicating 
that we use the ingredients coming from 
$F(w)$ and $F(w)^{\prime}$ in $\formula[{u,w}]^{\prime}$. Thus, 
$\formula[{u,w}]=\formula[{u,w}]^{\prime}$ which is what we wanted to show.

We next claim that $(u,w)$ stays the same under the moves 
in \autoref{Eq:ProofPhantom} and the consequent 
reorientation. Indeed, there is a bijection of 
the bases on either side since the counting of flows 
(another way to index $\mathrm{rank}_{\Zq}\big(W(u,w)\big)$) comes out 
to be the same. To see this, recall that a flow is a labeling of 
the edges of thickness $k$ with a $k$ element subset of $\{1,\dots,n\}$ 
such that locally {\eg}
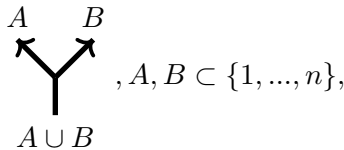
\begin{gather*}
\begin{tikzpicture}[anchorbase,scale=1,yscale=-1]
\draw[usual,directed=0.99] (0.5,0.5) to (0,0)node[above]{$A$};
\draw[usual,directed=0.99] (0.5,0.5) to (1,0)node[above]{$B$};
\draw[usual] (0.5,1)node[below]{$A\cup B$} to (0.5,0.5);
\end{tikzpicture}
,A,B\subset\{1,\dots,n\},
\end{gather*}
holds.
Then we first note that changing 
the orientation does not change the number of flows. 
The operations in \autoref{Eq:ProofPhantom} also do not change
the flow, so the claim follows.

Finally, we need to argue that the rank of the 
$\gln$ web algebra does not change when using phantom edges and crossings 
freely. To this end, note that \autoref{L:FFormPhantom} implies 
that webs seen as objects in three space do not change under the usage 
of phantom crossings. In particular, the evaluation formula that gives 
$(u,w)$ does not actually see them at all.

\textit{(Part b.)} By \autoref{T:EvaluationMain} we have that the
formula in \autoref{Eq:EvaluationMain} matches 
the evaluation on webs. The latter, by \cite[Theorem 5.1]{MuOhYa-webs} 
(strictly speaking we use different conventions but the arguments of \cite[Theorem 5.1]{MuOhYa-webs} still apply in our conventions), is 
known to give a link invariant. That this invariant is, up to conventions again, the Reshetikhin--Turaev polynomial of interest
follows then from \cite[Corollary 6.2.3]{CaKaMo-webs-skew-howe}.
\end{proof}

\begin{proof}[Proof of ``The dual canonical webs, {\ie} \autoref{P:EvaluationDualCanonical}'']
Directly from \autoref{T:EvaluationMain} and
\cite[Theorem 4.19]{Tu-gln-bases}.
\end{proof}



\end{document}